\documentclass[11pt]{article}
\usepackage{amssymb}
\usepackage{graphicx}

\usepackage{amsmath, amsthm, amssymb}
\usepackage{enumerate}
\usepackage{graphicx}
\usepackage{float}
\usepackage{hyperref}
\usepackage[margin=1in]{geometry}
\usepackage{comment}
\usepackage{color}
\numberwithin{equation}{section}

\begin{document}
\newtheorem{theorem}{Theorem}[section]
\newtheorem{proposition}[theorem]{Proposition}
\newtheorem{lemma}[theorem]{Lemma}
\newtheorem{cor}[theorem]{Corollary}
\newtheorem{definition}[theorem]{Definition}
\newtheorem{assumption}[theorem]{Assumption}
\newtheorem{conjecture}[theorem]{Conjecture}
\newtheorem{remark}[theorem]{Remark}

\title{On the Hardy-Littlewood-Sobolev type systems}

\author{Ze Cheng$^{a}$ \footnote{Partially supported by NSF DMS-1405175.},
 Genggeng Huang$^{b}$ \footnote{Partially supported by NSFC-11401376 and China Postdoctoral Science Foundation 2014M551391.},
 and Congming Li$^{b,a}$ \footnote{Partially supported by NSFC-11271166 and NSF DMS-1405175.}
 \\
$^a$ Department of Applied Mathematics, \\
University of Colorado Boulder, CO 80309, USA \\
$^b$ Department of Mathematics, INS and MOE-LSC, \\
Shanghai Jiao Tong University, Shanghai, China
}

\maketitle
\date{}

\begin{abstract}
In this paper, we study some qualitative properties of Hardy-Littlewood-Sobolev type systems. The HLS type systems are categorized into three cases: critical, supercritical and subcritical. The critical case, the well known original HLS system, corresponds to the Euler-Lagrange equations of the fundamental HLS inequality. In each case, 
we give a brief survey on some important results and useful methods. Some simplifications and extensions based on somewhat more direct and intuitive ideas are presented. Also, a few new qualitative properties are obtained and several open problems are raised for future research. %Another aspect of this paper is to simplify some of the related methods based on somewhat more direct and intuitive ideas.
\end{abstract}

\section{Introduction to the Hardy-Littlewood-Sobolev type systems}
The well known Hardy-Littlewood-Sobolev inequality (HLS) states:
\begin{align}
    \int_{R^n}\int_{R^n} \displaystyle{\frac{f(x)g(y)}{|x-y|^{n-\gamma}}} dx dy
    \leq
    C(n,s,\gamma) ||f||_r||g||_s
 \end{align}
where $0 <\gamma < n$, $1< s, r  < \infty$, $\frac{1}{r} +
\frac{1}{s} + \frac{n-\gamma}{n}= 2$, $f \in L^r(R^n)$ and $g \in
L^s(R^n)$.

Define an operator $T$ such that $ Tg(x):=\displaystyle{\int}_{R^n} \displaystyle{\frac{g(y)}{|x-y|^{n-\gamma}} }dy$, $\gamma \in (0,n)$, then the HLS inequality becomes:
\begin{align}
	 ||Tg||_{\frac{ns}{n-s\gamma}} \leq C(n,s,\gamma) ||g||_s, \text{ or }
	 ||Tg||_p \leq C(n,s,\gamma) ||g||_{\frac{np}{n+\gamma p}},
\end{align}
where $\frac{n}{n-\gamma}<p< \infty$, and $1<s< n/{\gamma}$.

The best constant $C = C(n,s,\gamma)$ is the maximum of:
\begin{equation}\label{optimizingProblem}
    J(f,g) =  \int_{R^n}\int_{R^n}\frac{f(x)g(y)}{|x-y|^{n-\gamma}} dx dy
\end{equation}
with constraints $ \|f\|_r = \|g\|_s =1$.
The above optimizing problem leads us to a system of integral equations on $f$ and $g$.
Let $u = c_1 f^{r-1}$, $v = c_2 g^{s-1}$,
 $p=\frac{1}{r-1}$,
 $q=\frac{1}{s-1}$, and  choose suitable constants $c_1$ and
 $c_2$, we arrive at the
following system of Euler-Lagrange equations for the HLS
inequality :

\begin{equation} \left \{
 \begin{array}{l}
      u(x) = \int_{R^{n}} \frac{ v^q(y)}{|x-y|^{n-\gamma}} dy\\
      v(x) = \int_{R^{n}} \frac{ u^p(y)}{|x-y|^{n-\gamma}} dy
\end{array}
\right. \label{sys}
\end{equation}
with $ u, v >0$, $\;u \in
L^{p+1}$, $\;v \in L^{q+1}$, $0 < p < \infty$, $\;0 < q < \infty$,
$\;\frac{1}{p+1} + \frac{1}{q+1} =  \frac{n-\gamma}{n}$.

For the same $p,q,n,\gamma$, the solutions of  (\ref{sys}) also are solutions of (the reverse is nontrivial and might be false \cite{CLO05}):
\begin{equation} \left \{
 \begin{array}{l}
      (-\Delta)^{\gamma/2} u = v^q , \;\; u>0, \mbox{ in } R^n,\\
      (-\Delta)^{\gamma/2} v = u^p, \;\; v>0, \mbox{ in } R^n.
\end{array}
\right. \label{pde}
\end{equation}
For $0 < p, q < \infty$, we call \eqref{sys} and \eqref{pde} \textbf{Hardy-Littlewood-Sobolev (HLS) type systems}. This class of systems include some well-known special cases, for instance,
if $ p=q=\frac{n+\gamma}{n-\gamma}$, and $u(x)=v(x)$, the above system
becomes:

\begin{equation}
(-\Delta)^{\gamma/2} u = u^{(n+\gamma)/(n-\gamma)},\;\; u>0, \mbox{ in } R^n.
\label{pde10}
\end{equation}
In particular, when $n \geq 3 $, and $\gamma = 2,$
\begin{equation}
 -\Delta u =  u^{ (n+2)/(n-2) }, \;\; u>0, \mbox{
in } R^n. \label{pde2}
\end{equation}
Equation (\ref{pde10}) is equivalent to:
\begin{equation}
    u(x) = \int_{R^{n}} \frac{u(y)^{\frac{n+\gamma}{n-\gamma} }}{|x - y|^{n - \gamma}}
      dy,\;\; u >0 \mbox{ in } R^n,
      \label{sys10}
\end{equation}
for which we refer readers to Chen, Li and Ou \cite{CLO06}, where \eqref{pde10} is defined for non-integer $\gamma$.

%three cases
The HLS type systems can be categorized into three cases, namely, critical case
$\frac{1}{p+1} + \frac{1}{q+1} =  \frac{n-\gamma}{n}$, subcritical case $\frac{1}{p+1} +\frac{1}{q+1} > \frac{{n-\gamma}}{n}$, supercritical case $\frac{1}{p+1} +\frac{1}{q+1} < \frac{{n-\gamma}}{n}$. The HLS type systems behave very differently in each case, hence different methods are needed to study them.

 \textbf{The critical} case of the HLS type system, also known as the HLS system, was studied in an elegant paper of E. Lieb \cite{Lieb83}. The existence of solution is proved (see \cite{Lieb83} and details in Section 3), so the remaining question is the classification/uniqueness of positive solutions.
We conjecture that all positive solutions of (\ref{sys}) or (\ref{pde}) are given by a single
family of solutions via translation and scaling.
In other words, up to some translations and scalings,
the positive solutions of (\ref{sys}) or (\ref{pde}) are unique. This kind of uniqueness is called essential uniqueness. This problem is solved in some special cases but is still open in general.

In \textbf{the supercritical case}\,
%we study the existence and asymptotic analysis of solution.
the existence of solution is also established (see Section 3). %We refer to a theorem by Liu, Guo and Zhang \cite{LGZ06} (2006) in Section 3. %Their result is stated below,
%\begin{theorem}[Liu, Guo and Zhang 2006]
%The system:
%\begin{equation} \left \{
% \begin{array}{l}
%      (-\Delta)^{k} u = v^q , \;\; u>0, \mbox{ in } R^n,\\
%      (-\Delta)^{k} v = u^p, \;\; v>0, \mbox{ in } R^n.
%\end{array}
%\right.
%\end{equation}
%admits a positive solution in the critical and super-critical cases
%$\frac{1}{p+1}+\frac{1}{q+1}\leq\frac{n-2k}{n}$
%for any integer $k$ with $2k<n$. In fact, given $pq>1$, the above system admits a positive radial solution
%if and only if $\frac{1}{p+1}+\frac{1}{q+1}\leq\frac{n-2k}{n}$.
%\end{theorem}
The existence proof is based on a relatively new method which combines shooting method with degree theory. We give an outline of the method to prove the existence for an even more general system in Section 4.1, which contains supercrtical and critical HLS type systems as special cases.

We are also interested in asymptotic analysis of the solutions of supercritical HLS type systems. For example, what are the asymptotic expansions of the radial solutions?
Are all radial solutions scaling related? %What is the limit solution when the exponents approaches the critical level?
%For the last problem, \textbf{ if we fix $1=u(0)\geq v(0)$ then the limit has finite energy and is radial}.
These questions are addressed in Section 4.2.

Last, we consider \textbf{the subcritical} HLS type systems, in particular, the Lane-Emden system.
The so-called Lane-Emden conjecture states that,
for $0<p,q<\infty$, $\frac{1}{p+1} +\frac{1}{q+1} > \frac{n-2}{n}$,
\begin{equation}
 \left \{
 \begin{array}{l}
      -\Delta u(x) = v^q(x) , \;\; u\geq0, \mbox{ in } R^n,\\
      -\Delta v(x) = u^p(x) , \;\; v\geq0, \mbox{ in } R^n,
\end{array}
\right. \label{LaneEmden}
\end{equation}
has $u=0$ and $v=0$ as the unique locally bounded solution.

The conjecture naturally generalize to the systems \eqref{sys} or \eqref{pde} in the subcritical cases with an additional condition,
$pq>1$. Notice that $pq>1$ is a necessary condition for this conjecture to hold in high order HLS type systems. For example if $p=q=1$ and $\gamma=4$ we have solution $u=v=e^{w\cdot x}$ to \eqref{pde} for $w\in\mathbb{R}^n$ with $|w|=1$.

This paper is organized as following. In Section 2 we provide some useful estimates. In Section 3, we discuss about the classification of solution to critical HLS type system. In Section 4, existence and asymptotic analysis of solutions to supercritical cases are investigated. In Section 5, we focus on subcritical cases, in particular, the Lane-Emden system. Throughout this paper, we assume that$0<p,q<\infty$ and $pq>1$, and all positive solution are assumed locally bounded, unless specified otherwise.

\section{Some basic estimates}
For $pq>1$, denote the scaling component of system \eqref{LaneEmden} by
	\begin{align}
	\alpha = \dfrac{2(p+1)}{pq-1}, \ \ \beta = \dfrac{2(q+1)}{pq-1}.
	\end{align}
	
	\begin{remark} For $pq>1$, the critical hyperbola has a new form in terms of $\alpha$ and $\beta$,
		\begin{align*}
			\frac{1}{p+1} +\frac{1}{q+1} =
			\frac{n-2}{n} \Leftrightarrow \alpha +\beta = n-2.
		\end{align*}
		Hence, supercritical condition $\frac{1}{p+1} +\frac{1}{q+1} <
					\frac{n-2}{n}\Leftrightarrow \alpha +\beta < n-2$, and subcritical condition $\frac{1}{p+1} +\frac{1}{q+1} >
										\frac{n-2}{n}\Leftrightarrow \alpha +\beta > n-2$.
	\end{remark}
Here we present some useful estimates including comparison principle and energy estimates for \eqref{LaneEmden}. These estimates are valid for all three cases, i.e., critical, supercritical and subcritical. They are useful in many aspects, such as in asymptotic analysis for solutions to critical and supercritical cases, and to prove Liouville type theorem in subcritical case (Actually, to prove Lane-Emden conjecture a proper energy estimate is the key, which is detailed in Section 5).
These estimates are:
\begin{lemma}\label{energyEstimates}
Let $p,q>0$ with $pq>1$. For any positive (locally bounded) solution $(u,v)$ of \eqref{LaneEmden}
\begin{align}\label{uvEstimate1}
	\int_{B_R} u \leq C R^{n-\alpha} \ \text{and } \int_{B_R} v \leq C R^{n-\beta},
\end{align}
\begin{align}\label{uvEsitmatePQ}
	\int_{B_R} u^q \leq C R^{n-q\alpha} \ \text{and } \int_{B_R} v^p \leq C R^{n-p\beta},
\end{align}
\end{lemma}
and by Maximum-principle we get,
\begin{lemma}[Comparison principle]\label{comparisonPrinciple}
	Let $p\geq q >0$ with $pq>1$. Let $(u,v)$ be a positive bounded solution of \eqref{LaneEmden}. Then we have the following comparison principle,
	\begin{align*}
		\frac{v^{p+1}}{p+1}\leq \frac{u^{q+1}}{q+1} , \  x\in \mathbb{R}^n.
	\end{align*}
\end{lemma}
\emph{Do similar estimates exist for general HLS type systems \eqref{sys} and \eqref{pde}?} This is also an interesting question, yet as far as we know there has not been any answer to it.

\noindent Proof of comparison principle Lemma \ref{comparisonPrinciple}.

Let $l=(\frac{p+1}{q+1})^{\frac{1}{p+1}}$, $\sigma=\frac{q+1}{p+1}$. So $l^{p+1}\sigma=1$, and $\sigma\leq 1$. Denote $\omega = v- lu^{\sigma}.$
Since
\begin{align*}
	\Delta\omega &= \Delta v - l \nabla\cdot(\sigma u^{\sigma-1}\nabla u)  \\
					&= \Delta v - l\sigma(\sigma-1)|\nabla u|^2 -l\sigma u^{\sigma-1}\Delta u\\
					&\geq -u^q+l\sigma u^{\sigma-1}v^p \\
					&= u^{\sigma-1}((\frac{v}{l})^p - u^{q+1-\sigma}) \\
					&= u^{\sigma-1}((\frac{v}{l})^p - u^{\sigma p}),		
\end{align*}
we see $\Delta\omega \geq Cw^p >0$ if $w>0$.
Now, suppose $w>0$ for some $x\in\mathbb{R}^n$.

Case 1: $\exists x_0\in \mathbb{R}^n$ such that $\omega(x_0)=\displaystyle \sup_{\mathbb{R}^n} \omega(x)>0$. Then we have
$\Delta\omega(x_0) \leq 0.$
This contradicts with the fact that $\Delta\omega(x_0) > 0$ since $\omega(x_0) > 0.$

Case 2:There exists a sequence $\{x_R\}$, such that $\displaystyle \lim_{R\rightarrow\infty} \omega(x_R) = \displaystyle \sup_{\mathbb{R}^n} \omega(x) >0$. %We claim $\Delta \omega(x_R) = O(\frac{1}{R^2}).$

Define $\omega_R(x)=\phi(\frac{x}{R})\omega(x)$ where $\phi(x)$ is a cutoff function on $\mathbb{R}^n$ such that $\phi(x)\equiv 1$ on the unit ball $B_1(0)$ and $\phi(x)\equiv 0$ outside $B_2(0)$.  Take $x_R$ such that $\omega_R(x_R) =  \max_{\mathbb{R}^n} \omega_R(x)$.

We see that $\lim_{R\rightarrow\infty} \omega_R(x_R)=\lim_{R\rightarrow\infty} \omega(x_R) = \displaystyle \sup_{\mathbb{R}^n} \omega(x) >0$ and:
\begin{align*}
	& 0=\nabla \omega_R(x_R) = \phi(\frac{x_R}{R})\nabla\omega(x_R) + \frac{1}{R}\nabla\phi(\frac{x_R}{R})\omega(x_R) \\
	\Rightarrow & \nabla\omega(x_R) = O(\frac{1}{R}), \ \text{as } R\rightarrow +\infty.
\end{align*}
So,
\begin{align*}
	0 \geq \Delta\omega_R(x_R)&=\frac{1}{R^2}\Delta\phi(\frac{x_R}{R})\omega(x_R)+\frac{2}{R}\nabla\phi(\frac{x_R}{R})\cdot\nabla\omega(x_R) \\
	& +\phi(\frac{x_R}{R})\Delta\omega(x_R) \\
	\Rightarrow  \Delta\omega(x_R) \leq O(\frac{1}{R^2})
\end{align*}
This contradicts with the relation $\Delta\omega \geq Cw^p >0$ if $w>0$ (suppose $\lim_{R\rightarrow\infty}w(x_R)>C_0>0$ then $\lim_{R\rightarrow\infty}\Delta\omega > C_1>0$ for some $C_1$).

So, we see that $\omega\leq 0$ on $\mathbb{R}^n$, and this finishes the proof. $\Box$

\quad

Lemma \ref{energyEstimates} is first obtained by Serrin and Zou \cite{SZ96} (1996). In \cite{CHL14} a simpler proof is given, to which we refer readers for detail. Here we only sketch the proof.
First, we multiply \eqref{LaneEmden} with $\phi$, the first eigenfunction of $-\Delta$ operator on $B_R$ with eigenvalue $\lambda\sim O(\frac{1}{R^2})$. Then integrate the equations by parts and use the fact that $\phi$ is positive in $B_R$ (which leads to $\frac{\partial \phi}{\partial n}<0$ on $\partial B_R$ by Hopf's Lemma), and we obtain the following inequalities,
\begin{align}\label{uqBoundedByV}
	\int_{B_R} \phi u^q \leq \lambda \int_{B_R} \phi v, \quad
	\int_{B_R} \phi v^p \leq \lambda \int_{B_R} \phi u.
\end{align}
Apply comparison principle lemma \ref{comparisonPrinciple} to the first inequality above to get
	\begin{align*}
		\frac{1}{R^2} \int_{B_R} \phi v\geq C \int_{B_R} \phi v^{\frac{q(p+1)}{q+1}}.
	\end{align*}
Notice that $\frac{q(p+1)}{q+1}>1$ as $pq>1$, so by H\"{o}lder inequality
\begin{align*}
	&\frac{1}{R^2} \int_{B_R} \phi v\geq C (\int_{B_R} \phi v)^{\frac{q(p+1)}{q+1}} R^{-n\frac{qp-1}{q+1}} \\
	\Rightarrow & \int_{B_R} \phi v \leq C R^{n-\beta}.
\end{align*}
Therefore, by \eqref{uqBoundedByV}
\begin{align*}
	\int_{B_R} \phi u^q \leq C R^{n-\beta-2} =CR^{n-q\alpha},
\end{align*}
and \textbf{Case 1: $q \geq 1$}, by H\"{o}lder's inequality:
\begin{align*}
	\int_{B_R} \phi u \leq (\int_{B_R} \phi u^q)^{\frac{1}{q}}(\int_{B_R} \phi )^{\frac{1}{q'}}
			\leq CR^{\frac{n}{q}-\alpha} R^{\frac{n}{q'}}
			= CR^{n-\alpha},
\end{align*}
and by the second inequality of \eqref{uqBoundedByV},
\begin{align*}
	\int_{B_R} \phi v^p \leq CR^{n-\alpha-2} = CR^{n-p\beta}.
\end{align*}
\textbf{Case 2: $q<1$}, this case is more complex than the first one. By the first equation of \eqref{LaneEmden} we have $-\Delta u\leq 0$, then multiply it with $\eta^2 u^{\gamma}$ where $\eta \in C^{\infty}_0(\mathbb{R}^n)$ and $\eta\in(0,1)$ and integrate over whole space, we get
\begin{align*}
	\int_{B_R} \frac{4}{\gamma^2}|D(u^{\frac{\gamma}{2}})|^2\leq \frac{C_\gamma}{R^2}\int_{B_{2R}}u^{\gamma}.
\end{align*}
Now we use Poincar\'{e} inequality to induce an embedding inequality on the support of $\eta$.
Then we use H\"older inequality and estimate for $\int_{B_R} u^q$ in previous proof to obtain,
\begin{align*}
	\int_{B_R} u^{\frac{n}{n-2}\gamma} \leq  C R^{n-\gamma\frac{n}{n-2}\alpha}.
\end{align*}
Then the proof is finished by taking $\gamma\leq q$ then use H\"older inequality to get $$\int_{B_R} u^\theta \leq CR^{n-\theta\alpha,}$$ for any $\theta\in (0,\frac{n}{n-2})$.

\section{Existence and classification of solutions for critical HLS type systems}
The existence of solution to critical HLS type systems, for both integral equations \eqref{sys} and PDE \eqref{pde}, is completely revolved.

In \cite{Lieb83}  E. Lieb (1983) established the existence of ground state solution (i.e. the optimizer of variational problem \eqref{optimizingProblem}). Later, people find that shooting method is a powerful tool to prove existence of solution to both critical and supercritical cases of \eqref{pde} with integer power of Laplacian. We sketch the proof of the shooting method in Section 4.1.

In the following system, Liu, Guo and Zhang \cite{LGZ06} (2006) have obtained the existence for both critical and supercritical cases,
\begin{theorem}\label{LGZ}
The system
\begin{equation} \left \{
 \begin{array}{l}
      (-\Delta)^{k} u = v^q , \;\; u>0, \mbox{ in } R^n,\\
      (-\Delta)^{k} v = u^p, \;\; v>0, \mbox{ in } R^n,
\end{array}
\right.
\end{equation} \label{HLSpde}
admits a positive solution in the critical and super-critical cases
$\frac{1}{p+1}+\frac{1}{q+1}\leq\frac{n-2k}{n}$
for any integer $k$ with $2k<n$. In fact, given $pq>1$, the above system admits a positive radial solution
if and only if $\frac{1}{p+1}+\frac{1}{q+1}\leq\frac{n-2k}{n}$.
\end{theorem}
In Section 4.1, we outline the proof of a more general theorem which contains the result above.
In the scalar case:
\begin{equation}
(-\Delta)^{k} u = u^p,\;\; u>0, \mbox{ in } R^n,
\label{pde10s}
\end{equation}
where $p>1$, $2k<n$. Lei and Li (2013) \cite{LL13} showed that
\begin{theorem}
For $p>1$, equation \eqref{pde10s} admits a locally bounded solution if and only if $p \geq \frac{n+2k}{n-2k}$.
\end{theorem}

\begin{remark}
	Notice that $p>1$ is necessary for the ``only if'' part to be true in the theorem above. If $p=1$ and $k$ is even, \eqref{pde10s} has many solutions. Does \eqref{pde10s} have solution when $k>1$ and $k$ is odd? This is an open question.
\end{remark}
%As for \eqref{sys}, the existence of ground state solution (i.e. the optimizer of variational problem \eqref{optimizingProblem}) in the critical case was established by E. Lieb  \cite{Lieb83} (1983). 
Moreover,
Lei and Li (2013) \cite{LL13} showed that
\begin{theorem} \label{th1.2}
Assume $pq>1$, then the HLS type system (\ref{sys})
has a pair of positive solutions $(u,v) \in L^{p+1}(R^n)
\times L^{q+1}(R^n)$ if and only if it is critical:
\begin{equation} \label{1.3ksce}
\frac{1}{p+1}+
\frac{1}{q+1}=\frac{n-\gamma}{n}.
\end{equation}

\end{theorem}

Now we shall focus on the classification of the solutions. As mentioned in the introduction, for critical HLS type systems, we conjecture that\textbf{ all positive solutions of (\ref{sys}) or (\ref{pde}) are given by a single
family of solutions via translation and scaling.}
As a special case of this conjecture, E. Lieb raised an open problem in \cite{Lieb83} (1983) that asks for $p=q$ if \eqref{sys10} has unique solution up to scaling and translating. Lieb's problem was completely resolved with the introduction of the integral form of method of moving planes by Chen, Li and Ou \cite{CLO06} (2006). The conjecture is still open in general.

The method of moving planes was introduced by A.D. Alexandrov in 1950s, and then developed by J. Serrin \cite{Serrin71} (1971) and Gidas, Ni, Nirenberg \cite{GNN81} (1981). Caffarelli, Gidas and Spruck classified all the solutions to \eqref{pde2} in \cite{CGS89} (1989) by the method of moving planes. Then Chen and Li simplified their proof \cite{CL91} (1991). Wei and Xu generalized this result to higher order conformally invariant equations \cite{WX99} (1999). Then Chen, Li and Ou finally implemented the method of moving planes to the integral equation \eqref{sys} with $p=q$ and solved Lieb's problem \cite{CLO06} (2006).

%method of MP no integrable solution

Hereinafter, by ground state solutions we mean the solution $(u,v)$ corresponding to the optimizer $(f,g)$ of the functional \eqref{optimizingProblem}. To tackle the problem, we break it down to the following open problems:
\begin{enumerate}
\item \emph{all ground state solutions are translation and scaling related;
\item  all finite energy solutions are ground state;
\item  all radial solutions are finite energy solution (the converse is known to be true, see \cite{CLO06});
\item  all solutions decaying to zero at infinity are radial;
\item   all bounded solutions are radial.}
\end{enumerate}

%By Pohozaev identity, Serrin and Zou, Mitidieri, Lei and Li

\section{The super-critical HLS type systems}
In this section, we study the existence and asymptotic analysis of solution to supercritical HLS type systems.

\subsection{Existence of solution to critical and supercritical HLS type systems}

% % % % % % % % % % % % % % % % % % % % % % % Need to revise

For the Lane-Emden system \eqref{LaneEmden}, Serrin and Zou (1998) \cite{SZ98} used shooting method to obtain the existence of solution.
Liu, Guo and Zhang (2006) \cite{LGZ06} introduced a degree approach to shooting method (see also Li (2011) \cite{LiarXiv13}) to obtain radial positive (locally bounded) solution for both critical and supercritical cases of HLS type system \eqref{pde} in a uniform way.

 In short, the degree approach of shooting method combines three ingredients together: degree theory, target map (i.e. shooting) and non-existence on balls with Dirichlet boundary condition (Pohoza\'ev identities). Approach of this kind, relating the existence of solutions in $\mathbb{R}^n$ to the non-existence to a corresponding Dirichlet problem on balls, is implemented by many mathematicians, for instance earlier by Berestycki, Lions and Peletier \cite{BLP81} (1981).

% % % % % % % % % % % % % % % % % % % % % % % % %

Here we present a theorem by Cheng and Li (2015) \cite{CL15}, which implements degree approach of shooting method to Sch\"odinger type systems with sign-changing nonlinearities (HLS type system \eqref{pde}  is included as a special case). Consider a very general system in whole space,
\begin{align}\label{original}
  \left\{\begin{array}{cl}
      	-\Delta u_{i} = f_{i}(u) & \text{in } \mathbb{R}^n, \\
	       u_{i} > 0             & \text{in } \mathbb{R}^n, \\
 	     % u_{i}(x) \rightarrow 0 & \text{uniformly as } |x| \rightarrow \infty
        \end{array}
 \right.
\end{align}
and its corresponding local Dirichlet problem,
\begin{align}\label{Dirichlet}
  \left\{\begin{array}{cl}
      	-\Delta u_{i} = f_{i}(u) & \text{in } B_R, \\
	       u_{i} > 0             & \text{in } B_R, \\
           u_i=0                 & \text{on } \partial B_R,
 	     % u_{i}(x) \rightarrow 0 & \text{uniformly as } |x| \rightarrow \infty
        \end{array}
 \right.
\end{align}
where $B_R=B_R(0)$ for any $R>0$ and $i = 1,2,\cdots, L$. We have
\begin{theorem}\label{existencetheorem}
Given the nonexistence of radial solution to system \eqref{Dirichlet} for all $R>0$, the system \eqref{original} admits a radially symmetric solution of class $C^{2,\alpha}(\mathbb{R}^{n})$ with $0<\alpha<1$, if $f= (f_1(u), f_2(u),\cdots, f_L(u)):\mathbb{R}^L\rightarrow \mathbb{R}^L$ satisfies the following assumptions:
\begin{enumerate}
  %\item $f=\nabla F$ for some $F(u)\geq 0$ in $\mathbb{R}^{L}$, and $F$ is homogeneous;
  \item $f$ is continuous in $\overline{\mathbb{R}^{L}_{+}}$
and locally Lipschitz continuous in $\mathbb{R}^{L}_{+}$, and furthermore,
  \begin{align}\label{decayAssumption}
        \sum_{i=1}^L f_i(u) \geq 0  \text{ in } \mathbb{R}_+^L;
  \end{align}
  \item If $\alpha \in \partial\mathbb{R}^{L}_{+}$ and $\alpha\neq 0$, i.e., for some permutation $(i_1,\cdots,i_L)$, $\alpha_{i_1}=\cdots=\alpha_{i_m}=0$, $\alpha_{i_{m+1}},\cdots, \alpha_{i_L}>0$ where $m$ is an integer in $(0,L)$, then $\exists \delta_0 = \delta_0(\alpha) > 0$ such that for $\beta \in \mathbb{R}^{L}_{+}$ and $|\beta-\alpha|<\delta_0$,
  \begin{align}\label{ControlInequality}
    \sum_{j=m+1}^L |f_{i_j}(\beta)| \leq C(\alpha) \sum_{j=1}^m f_{i_j}(\beta),
  \end{align}
  where $C$ is a non-negative constant that depends only on $\alpha$.
\end{enumerate}
\end{theorem}
Notice that the nonlinear terms of \eqref{HLSpde} (i.e. the HLS type system \eqref{pde} with $\gamma=2k$) satisfy the assumptions \eqref{decayAssumption}-\eqref{ControlInequality}. So Theorem \ref{existencetheorem} recovers Theorem \ref{LGZ} given the non-existence of solution to the corresponding Dirichlet problems. For nonexistence part, people usually apply Pohoza\'ev inequalities \cite{Pohozaev65} (see also \cite{PS86}). This is now more or less standard, a large amount of literature can be found in related topics, for example \cite{Mitidieri93} and \cite{QS07}.

In what follows we outline the proof of Theorem \ref{existencetheorem}, which should exhibit the standard procedure of the degree approach of shooting method.

In the view of seeking radial solutions of \eqref{original}
we solve the following initial value
problem with any initial value $\alpha=(\alpha_1, \cdots, \alpha_L)$ with $\alpha_i>0, i=1,2,....L.$ Denote the solution as $u(r,\alpha)$:
\begin{equation} \label{eq:ode}
\left\{ \begin{aligned}
        &u_i^{''}(r) +\frac{n-1}{r} u_i^{'}(r)= -f_i(u) \\
         &u_i'(0)=0, u_i(0)=\alpha_i \quad i=1,2,\ldots,L.
                          \end{aligned} \right.
\end{equation}
We want to find the suitable initial value $\alpha$ so that $u_i(r,\alpha) >0$ for all $r>0$.

When $L=1$ the question is simple, the assumption that \eqref{Dirichlet} admits no solution
is equivalent to $u_1(r,\alpha) >0$ for all $r>0$. Then there exists a global solution for any initial value, and we are done.

When $L \geq 2$, instead of one dimensional initial value which scales to each other,
we are encountered with multi-dimensional initial value. Among $\alpha_i$'s,
in many critical cases as well as in many supercritical cases,
there is at most one scaling class (one-dimensional) of initial values
from which we can shoot to a global solution. To show the existence of positive solutions of
(\ref{eq:ode}), up to a simple scaling, we have to find the special 1-D initial values.
This is the main reason why there are so many results  in the scalar case
but very little about \eqref{original} for a long time period.

The degree theory approach for the shooting method gives a simple solution to this difficult problem.
It can be used to solve a much larger class of problems.
The argument starts with
defining the target map $\psi$. By $\alpha>0$ we mean $\alpha$ is an interior point of $\mathbb{R}^L_+$, and let $r_0$ be the smallest value of $r$ for which $u_i(r, \alpha)=0$ for some $i$. %If $r_0=\infty$ then we are done. So, we assume that $r_0<\infty$.
Define a map
\begin{align*}
	\psi: & R^L_+ \rightarrow \partial R^L_+	 \\
		& \alpha \mapsto \psi(\alpha)=\left\lbrace \begin{array}{ll}
		u(r_0, \alpha), &\text{ if } \alpha > 0, r_0<\infty,\\
		\displaystyle\lim_{r \rightarrow \infty} u(r, \alpha), &\text{ if } \alpha > 0, r_0=\infty, \\
		\alpha, &\text{ if } \alpha\in \partial R_+^L.
		\end{array}
		\right.
\end{align*}

%This and the assumption that $F(\beta) \neq 0$ when $\beta\gvertneqq 0$ ensure that $\psi(\alpha) \in \partial R_+^L$.

Then we need to show that $\psi$ is continuous from $R_+^L$
to $ \partial R_+^L$.
Assumptions \eqref{decayAssumption}-\eqref{ControlInequality} guarantee this.

In the next step, applying the degree
theory, we show that $\psi$ is onto from $A_a$ to $B_a$ where:

\begin{equation} \label{eq:scalar1}
\left\{ \begin{aligned}
 & A_a\triangleq \{ \alpha \in R_+^L \,\, {\bf \mid } \quad {\tiny \displaystyle\sum_{i=1,\cdots, L}} \alpha_i=a\},\\
  & B_a\triangleq\{ \alpha \in \partial R_+^L \,\, {\bf \mid } \quad {\tiny \displaystyle\sum_{i=1,\cdots, L}} \alpha_i \leq a\},
                          \end{aligned} \right.
\end{equation}
for any $a>0$.
In particular, there exists at least one
$\alpha_a \in A_a$ for
sufficiently small $a>0$ such that $\psi(\alpha_a)=0$.

Shooting from the initial value $\alpha_a$, by the assumption that the system
\eqref{Dirichlet} admits no radially symmetric solution, we obtain a solution of \eqref{original} which decays to 0 at infinity ($\displaystyle\lim_{r \rightarrow \infty} u(r, \alpha_a)=\psi(\alpha_a)=0$). Notice that we did not rule out the radial solution that does not decay at infinity for system with sign-changing nonlinearity, however for positive nonlinear source term such as HLS type system \eqref{pde}, the radial solution must monotone decrease, hence its radial solution must decay to zero at infinity.

We remark that assumption \eqref{ControlInequality} is a necessary condition for the target map to be continuous. A trivial degenerate case in the following does not satisfy \eqref{ControlInequality},
\[ \left\{ \begin{array}{l}
         -\triangle u_1 = u_1^p, \\
                  -\triangle u_2=u_2^q,
                  \end{array} \right.
         \]
which can be decoupled to the study of two scalar equations.
If $p> \frac{n+2}{n-2}$ and $q< \frac{n+2}{n-2}$ then the target map is not continuous at
$(0,a)$ for any $a>0$. In fact, $\phi(0,a)=(0,a)$ and $\psi(\delta, a)=(h(\delta),0)$ where $0<h(\delta)<\delta$
and $\psi(\delta, a)\rightarrow (0,0) \neq (0,a)$ as $ \delta \rightarrow 0^+$.

\subsection{Asymptotic analysis}

%J. Lim and C. Li, C. Jin, W. Chen, Chen-Li-05 best integrability
%J. Lim and li critical asymptotic

% % % % % % % % % % % % %Question

In the asymptotic analysis we try to answer two questions: What are the asymptotic expansions of the radial solutions? Are all radial solutions scaling related?

%\begin{itemize}
%  \item What are the asymptotic expansions of the radial solutions?
%  \item  Are all radial solutions scaling related?
%  \item  What is the limit solution when the exponents approaches the critical level?
%\end{itemize}

%For the last problem,  if we fix $1=u(0)\geq v(0)$ then the limit has finite energy and is radial.

% % % % % % % % % %Lei Li Ma

In \cite{LLM12}, Lei, Li and Ma (2012) considered a more general system as following,
\begin{equation}
 \left \{
   \begin{array}{l}
      u(x) = \frac{1}{|x|^{\alpha}}\int_{R^{n}} \frac{v(y)^q}{|y|^{\beta}|x-y|^{\lambda}}  dy\\
      v(x) = \frac{1}{|x|^{\beta}}\int_{R^{n}} \frac{u(y)^p}{|y|^{\alpha}|x-y|^{\lambda}}  dy
   \end{array}
   \right.\label{ELWHLS}
 \end{equation} where $\alpha+\beta+\lambda \leq n$, and
     \begin{equation}
     \left \{
   \begin{array}{l}
   u, v \geq 0, \;\; 0 <p, q<\infty,\;
    0 <\lambda < n, \; \alpha+\beta  \geq 0,\; \\ % \label{ELWHLScond1}
   \frac{\alpha}{n}  <
  \frac{1}{p+1}<\frac{\lambda+\alpha}{n},\;\frac{1}{p+1}+\frac{1}{q+1}=\frac{\lambda+\alpha+\beta}{n}.
\end{array}
   \right.
 \label{ELWHLScond3}
   \end{equation}
Notice that HLS type system \eqref{sys} is a special case ($\alpha=\beta=0$) of \eqref{ELWHLS}.
%\begin{definition}
%A function $u$ is asymptotic to $\frac
%{A}{|x|^s}$ near $x=0$, that is
%\begin{equation}
%u(x) \simeq \frac {A}{|x|^s}\;{\rm at}\; |x| \simeq 0, \quad {\rm
%if} \lim_{|x|\to 0}|x|^s u(x) = A, \label{def_1}
%\end{equation}
%for a positive number $s$ and $A \neq 0 \;(\infty)$.
%\end{definition}

%In a similar way, we define near $x=\infty$ that,
\begin{definition}
 A function $u$ is asymptotic to $\frac {B}{|x|^t}$ near
$x=\infty$, that is
\begin{equation}
u(x) \simeq \frac {B}{|x|^t}\;{\rm at} \;|x| \simeq \infty, \quad
{\rm if} \lim_{|x|\to \infty}|x|^t u(x) = B,\label{def_2}
\end{equation} for a positive number $t$ and $B \neq 0\; (\infty)$.
\end{definition}

Lei, Li and Ma proved that
%\begin{theorem}
%\it Let $(u,v) \in L^{p+1}(R^n)\times
%L^{q+1}(R^n)$ be a pair of positive solutions of the system
%(\ref{ELWHLS}) with (\ref{ELWHLScond3}). Suppose that $p \geq 1,\;
%q \geq 1,\; pq \neq 1$ and $\alpha+ \beta \geq 0$. If
%$\lambda+(q+1)\beta<n$, then for small $|x|$ we have
%\begin{equation}
%   u(x) \simeq \frac {A_0}{|x|^\alpha}, \quad \label{uasymz}
% \end{equation}
%and
%\begin{eqnarray}
%   v(x)& \simeq &
%       \left\{\begin{array}{cll}
%   \displaystyle\frac {A_1}{|x|^\beta}, \quad
%       &{if} \;\lambda +\alpha(p+1)<n \\
%         \displaystyle\frac {A_2 |\ln |x||}{|x|^\beta}, \quad
%         &{if}\; \lambda+\alpha(p+1)=n \\
%         \displaystyle\frac {A_3} {|x|^{ \alpha(p+1)+\beta+\lambda-n}}, \quad
%         &{ if}\; \lambda +\alpha(p+1)>n  \label{vasymz}
%    \end{array}\right.
%    \end{eqnarray}
%where $A_0 =\int_{R^n}\frac{v^q(y)} {|y|^{\lambda+\beta}} \,dy, $
%$A_1 =\int_{R^n} \frac{u^p(y)}{|y|^{\lambda+\alpha}} \, dy,\;$
%   $A_2 =  |S^{n-1}|(\int_{R^n}\frac{v^q(y)} {|y|^{\lambda+\beta}} \,dy)^p$ \\
%   and
%   $ A_3 =(\int_{R^n}\frac{v^q(y)} {|y|^{\lambda+\beta}} \,dy)^p
% \int_{R^n}\frac{dz}{|z|^{\alpha(p+1)}|e-z|^{\lambda}}$
%($e$ is a unit vector in $R^n$ and $|S^{n-1}|$ is the surface
%area of the unit sphere).
%
%\end{theorem}

\begin{theorem}
\it Let $(u,v) \in L^{p+1}(R^n)\times
L^{q+1}(R^n)$ be a pair of positive solutions of the system
(\ref{ELWHLS}) with (\ref{ELWHLScond3}). Suppose that $p \geq 1,\;
q \geq 1,\; pq \neq 1$ and $\alpha+\beta \geq 0$. If $\lambda
q+\beta (q+1)>n$, then  we have
\begin{equation}
     u(x) \simeq  \frac {B_0}{|x|^{\lambda+\alpha}}, \text{ for large $|x|$,}\\
     \label{uasyminf}
\end{equation}

and

\begin{eqnarray}
     v(x) & \simeq & \left\{\begin{array}{cll}
          \displaystyle\frac {B_1}{|x|^{\lambda+\beta}},
              \quad   &{ if}\; \lambda p +\alpha(p+1)>n\\
        \displaystyle\frac {B_2 |\ln |x||}{|x|^{\lambda+\beta}},
             \quad &{ if}\; \lambda p +\alpha(p+1)=n \\
          \displaystyle\frac {B_3} {|x|^{(\alpha +\lambda)(p+1)+\beta-n}},
              \quad &  { if}\; \lambda p +\alpha(p+1)<n   \label{vasyminf}
    \end{array}\right. \text{ for large $|x|$}
 \end{eqnarray}
where $B_0 =\int_{R^n}\frac{v^q(y)} {|y|^{\beta}} \,dy$, $B_1
=\int_{R^n}\frac{u^p(y)} {|y|^{\alpha}} \,dy$, $B_2 =
|S^{n-1}|(\int_{R^n}\frac{v^q(y)} {|y|^{\beta}} \,dy)^p$ and\\
$B_3 = (\int_{R^n}\frac{v^q(y)} {|y|^{\beta}} \,dy)^p
\int_{R^n}\frac{dz}{|z|^{2n-(\alpha+\lambda)(p+1)}|e-z|^{\lambda}}$.

\end{theorem}
\par Recently, the authors consider the supercritical Lane-Emden system \eqref{LaneEmden}
%\begin{equation}\label{4.2.1}
%\begin{cases}
%\Delta u+ v^p=0,\\
%\Delta v+ u^q=0,
%\end{cases}\quad\text{in}\quad \mathbb R^n,
%\end{equation}
with $p,q$ satisfying the supercritical condition, i.e.
\begin{equation}\label{4.2.2}
\frac 1{p+1}+\frac 1{q+1}<1-\frac2{n}, \quad 0<p,q<+\infty.
\end{equation}
The existence of the positive radial symmetric solutions are already known by \cite{LGZ06} via shooting method. We consider the asymptotic property of these radial solutions in the supercritical cases.
\begin{theorem}
\label{thm4.2.1}
Suppose $(u,v)$ be locally bounded non-negative radial solutions of \eqref{LaneEmden} with $p,q$ satisfying \eqref{4.2.2}. Then either
\begin{equation}\label{4.2.13}
C_1 r^{-\frac{2(p+1)}{pq-1}}\leq u(r)\leq C_2 r^{-\frac{2(p+1)}{pq-1}},\quad C_1 r^{-\frac{2(q+1)}{pq-1}}\leq v(r) \leq C_2 r^{-\frac{2(q+1)}{pq-1}}
\end{equation}
for some positive constants $C_1,C_2>0$ and $r$ large or $u=v\equiv 0$.
\end{theorem}
A more interesting open problem is if $u(r) r^{\frac{2(p+1)}{pq-1}}, v(r)r^{\frac{2(q+1)}{pq-1}}$ converge as $r\rightarrow\infty$.
%Once it is confirmed, another question would be we can show the radial solutions of Lane-Emden system are essentially unique. 

Proof.
Recall that we define $\alpha=\frac{2(p+1)}{pq-1},\beta=\frac{2(q+1)}{pq-1}$ (see Section 2). %where
%$$\frac{1}{p+1}+\frac{1}{q+1}<\frac{n-2}{n}, \quad 0<p,q<+\infty.$$
Considering the radial solutions of \eqref{LaneEmden}, we have the following ODE system,
\begin{align}\label{odeLaneEmden}
\left\lbrace \begin{array}{cc}
 -(r^{n-1}u')'=r^{n-1}v^p, \\
 -(r^{n-1}v')'=r^{n-1}u^q,   
\end{array}\right. 
\end{align}
Since we consider the radial solutions, by Lemma \ref{energyEstimates}, the upper bound of \eqref{4.2.13} follows easily.

%Recall that for \eqref{LaneEmden} we have the following Pohozaev's identity
%\begin{equation}
%\begin{split}
%&(\frac n{p+1}-a_1)\int_{B_R} v^{p+1}+(\frac{n}{q+1}-a_2)\int_{B_R} u^{q+1}\\=
%&|S^{n-1}|R^{n-1}(\frac{Rv^{p+1}(R)}{p+1}+\frac{Ru^{q+1}(R)}{q+1}+Ru'v'(R)+a_1u'v(R)+a_2v'u(R))
%\end{split}
%\end{equation}\label{pohozaevId1}
%where $a_1+a_2=n-2$. By the supercritical condition, we have $\alpha+\beta<n-2$.
%By supercritical property, we have $a_1+a_2=n-2>\frac n{p+1}+\frac n{q+1}$. Without loss of generality, we may assume $p\geq q$. Then we can choose $a_1<a_2$ such that
%\begin{equation}
%\begin{split}
%a_2(uv)'(R)&\leq (a_2-a_1) u'v(R)-R[\frac{v^{p+1}(R)}{p+1}+\frac{u^{q+1}(R)}{q+1}+u'v'(R)]\\
%&\leq -R[\frac{v^{p+1}(R)}{p+1}+\frac{u^{q+1}(R)}{q+1}]\\
%&\leq -CR(uv)^\gamma, \quad \gamma=\frac{\alpha+\beta+2}{\alpha+\beta}>1,
%\end{split}
%\end{equation}
%where the last inequality above is obtained by Young's inequality. This implies
%$$[(uv)^{1-\gamma}]'\geq CR.$$
%Directly integrating from $0$ to $R$, we get $uv\leq CR^{-\alpha-\beta}$. By comparison principle (Lemma \ref{comparisonPrinciple}), we have $v\leq CR^{-\beta}$.
%Now we turn back to the ODE, we get
%\begin{equation}
%\begin{split}
%-r^{n-1}u'(r)=\int_0^r s^{n-1}v^{p}(s)ds\leq Cr^{n-\alpha-2}.
%\end{split}
%\end{equation}
%Then we have $-u'(r)\leq Cr^{-\alpha-1}$, and integrate from $R$ to $\infty$ to get $u(r)\leq Cr^{-\alpha}$.

%It remains to prove that $uv\geq Cr^{-\alpha-\beta}$.
%\par Step 2: 
\par It remains to prove that $u,v$ are positive solutions, then $\displaystyle\liminf_{R\rightarrow +\infty} R^\alpha u(R)=c>0$.
\par Suppose not, we may assume $u(r)=w(r)r^{-\alpha}$ and $\displaystyle\liminf_{r\rightarrow +\infty} w(r)=0$. Now, we derive the upper bound of $-u'(r), -v'(r)$. Again, by comparison principle, we have
$$v(r)\leq \sigma u^{\frac{q+1}{p+1}}=\sigma w^{\frac{q+1}{p+1}}(r)r^{-\beta},$$
where $\sigma=\sqrt[p+1]{\frac{p+1}{q+1}}$. As discussed in \cite{Mitidieri92}, we have
\begin{equation}
\begin{split}
-ru'(r)=(N-2)u(r)-(N-2)u(0)-\int_r^\infty sv^p(s)ds\leq (N-2)u(r)=(N-2)w(r)r^{-\alpha}
\end{split}
\end{equation}
The same conclusion also holds for $-rv'(r)$ such that $-rv'(r)\leq (N-2)\sigma w^{\frac{q+1}{p+1}}(r) r^{-\beta}$. 
Recall the following Pohozaev's identity for radial solutions,
\begin{equation}
\begin{split}
&(\frac n{p+1}-a_1)\int_0^R s^{n-1} v^{p+1}(s)ds+(\frac{n}{q+1}-a_2)\int_0^R s^{n-1}u^{q+1}(s)ds\\=
&R^{n-1}(\frac{Rv^{p+1}(R)}{p+1}+\frac{Ru^{q+1}(R)}{q+1}+Ru'v'(R)+a_1u'v(R)+a_2v'u(R)).
\end{split}
\end{equation}
We distinguish two cases:
\par 1. There exists $r_0>0$ such that for $r>r_0$, $w'(r)\leq 0$.
\par 2. For any $l>0$, there exists $s_l>l$ such that $w'(s_l)>0$. Since $\displaystyle\liminf_{r\rightarrow \infty} w(r)=0$, % by mean value theorem, there exists $t_l>l$ such that $w'(t_l)\leq 0$. So, 
we can get a sequence of local minimum points $r_l$ of $w$, such that $w(r_l)\rightarrow 0$ as $r_l\rightarrow +\infty.$
\par For both cases, we have
\begin{equation}\label{lowerboundAsymptotic}
-u'(r_l)=-w'(r_l)r^{-\alpha}_l+\alpha w(r_l)r_l^{-\alpha-1}\geq \alpha w(r_l)r_l^{-\alpha-1}.
\end{equation}
Now consider Pohozaev's identity at $r_l$ and $a_1=\frac n{p+1}$. By the supercritical property of $p,q$, we have
\begin{equation}
\begin{split}
\int_0^{r_l} r^{n-1}u^{q+1}(r)dr \leq C(w^{q+1}(r_l)+w^{1+\frac{q+1}{p+1}}(r_l))r_l^{n-2-\alpha-\beta}\leq Cw^{1+\frac{q+1}{p+1}}(r_l)r_l^{n-2-\alpha-\beta}
\end{split}
\end{equation}
where we have used $\frac{q+1}{p+1}<q$.  By comparison principle, we also have 
$$\int_0^{r_l} r^{n-1}v^{p+1}(r)dr \leq Cw^{1+\frac{q+1}{p+1}}(r_l)r_l^{n-2-\alpha-\beta}.$$
 Now we return to the first equation of \eqref{odeLaneEmden},
\begin{equation}\label{upperboundAsymptotic1}
\begin{split}
-r_l^{n-1}u'(r_l)&=\int_0^{r_l}s^{n-1} v^p(s)ds\\ 
&\leq \left(\int_0^{r_l} s^{n-1}v^{p+1}\right)^{\frac p{p+1}}\left(\int_0^{r_l} s^{n-1}ds\right)^{\frac 1{p+1}}\\ 
&\leq C w^{\frac p{p+1}+\frac{p(q+1)}{(p+1)^2}}(r_l)r_l^{n-2-\alpha}
\end{split}
\end{equation}
Combine \eqref{lowerboundAsymptotic} and \eqref{upperboundAsymptotic1} we get $$C w^{\frac p{p+1}+\frac{p(q+1)}{(p+1)^2}}(r_l)r_l^{n-2-\alpha} \geq -r_l^{n-1}u'(r_l)\geq \alpha w(r_l) r_l^{n-2-\alpha}.$$
From $w(r_l)\neq 0$, we obtain
\begin{equation}
w^{\frac{pq-1}{(p+1)^2}}(r_l)\geq \frac 1C
\end{equation}
This contradicts with $w(r_l)\rightarrow 0$.
Hence $u(r)\sim r^{-\alpha}$ near $\infty$. 
Then integrate the second equation of \eqref{odeLaneEmden} from 0 to $r$,
\begin{equation}
-r^{n-1} v'(r)=\int_0^r s^{n-1} u^q ds\geq Cr^{n-q\alpha}.
\end{equation}
Dividing the above inequality by $r^{n-1}$ and integrating from $r$ to $\infty$, we derive $v(r)\geq c r^{-\beta}$.
\par This ends the proof of the asymptotics of $u,v$ at $\infty$.
%\par Consider $\tilde u(t)=e^{\alpha t}u(e^t)$, $\tilde v(t)=e^{\beta t}v(e^t)$. Then 
%\begin{equation}
%\begin{cases}
%\tilde u''+(n-2-2\alpha)\tilde u'-\alpha(n-2-\alpha)\tilde u+\tilde v^p=0,\\
%\tilde v''+(n-2-2\beta)\tilde v'-\beta(n-2-\beta)\tilde v+\tilde u^q=0,
%\end{cases}\quad t\in (-\infty,+\infty)
%\end{equation}
%If we assume $\tilde u(t)\rightarrow C_1$ as $t\rightarrow +\infty$, then set $\tilde u_m(t)=\tilde u(t+t_m)$ for some $t_m\rightarrow +\infty$. Since $\tilde v_n,\tilde u_n$ are uniformly bounded, by elliptic estimates, we get $\tilde u_n,\tilde v_n$ locally uniformly converges in $C^{1,\gamma}$ by Sobolev embedding theorem. Again, by Schauder's estimates, $\tilde u_n$ converges to $C_1$ in $C^{2,\gamma}$ which implies $\tilde v_n$ converges to $C_2$.
$\Box$

\section{Subcritical HLS type systems: Liouville type theorem and the Lane-Emden conjecture}

%\subsection{The Lane-Emden conjecture}

For subcritical HLS type systems, various Liouville type theorems are obtained. This kind of results are often based on the study of the Lane-Emden system \eqref{LaneEmden}.
The Lane-Emden conjecture has been lasting unsolved for decades. Many mathematicians have contributed in this question, for example, the pioneer job done by Mitidieri (1992) \cite{Mitidieri92} (see also \cite{Mitidieri96}) which solves the Lane-Emden conjecture in radial case. For expository reference about the Lane-Emden conjecture, readers can check \cite{CHL14, Souplet09} and reference therein. Among these mathematical works we mention a couple of results below.

	For $pq\leq 1$, it is known that system \eqref{LaneEmden} has no positive classical supersolutions (see Serrin and Zou (1996) \cite{SZ96}).

For $n=3$,  the conjecture is solved by two papers. First, Serrin and Zou (1996) \cite{SZ96} proved that there is no positive solution with polynomial growth at infinity.
\begin{theorem}[Serrin-Zou-1996]\label{SZ06}
		Let $n=3$. Lane-Emden system \eqref{LaneEmden} admits no solution given the solution has at most polynomial growth at infinity.
	\end{theorem}
Then Pol\'{a}\v{c}ik, Quittner and Souplet (2007) \cite{PQS07} removed the growth condition. In fact, they proved that no bounded positive solution implies no positive solution.
\begin{theorem}[Pol\'{a}\v{c}ik-Quittner-Souplet-2007 ]\label{PQS09}
		Let $pq>1$. Assume that \eqref{LaneEmden} does not admit
		any bounded nontrivial (nonnegative) solution in $\mathbb{R}^n$.
		Let $\Omega\neq\mathbb{R}^n$ be a domain. There exists $C=C(n,p,q)>0$
		such that any solution $(u,v)$ of \eqref{LaneEmden} in $\Omega$
		satisfies
		\begin{align*}
			u(x)\leq C dist^{-\alpha}(x,\partial\Omega), \ x\in\Omega,
		\end{align*}
		and
		\begin{align*}
			v(x)\leq C dist^{-\beta}(x,\partial\Omega), \ x\in\Omega.
		\end{align*}
	\end{theorem}
	
	\begin{remark}
	In \cite{PQS07}, $p,q$ were assumed to be both $>1$, however, their proof is valid for $pq>1$ and can be verified directly.
	\end{remark}
	
%	\begin{remark}
%	Theorem \ref{PQS09} implies that
%	\begin{align*}
%		\text{Existence of solutioin. } \Rightarrow \text{Existence of bounded solution.}
%	\end{align*}
%	\end{remark}	
	
This result has two important consequences. One is that combining with Serrin and Zou's result, one can prove the conjecture for $n=3$.
	\begin{cor}
	Lane-Emden conjecture is true when $n=3$.
	\end{cor}
The other consequence is that one can reduce the Lane-Emden conjecture to the nonexistence of bounded positive solution.
Thus, hereinafter we always assume that $(u,v)$ are bounded.

For $n=4$, the conjecture is recently solved by Souplet (2009) \cite{Souplet09}. In \cite{SZ96}, Serrin and Zou used the integral estimates to derive the nonexistence results. Souplet further developed the approach of integral estimates and solved the conjecture for $n=4$ along the case $n=3$. In higher dimensions, this approach provides a new subregion where the conjecture holds, but the problem of full range in high dimensional space still seems stubborn.
\begin{theorem}[Souplet-2009]
	For subcritical Lane-Emden system \eqref{LaneEmden}, i.e. $p,q>0$ and $\frac{1}{p+1} +\frac{1}{q+1} > \frac{n-2}{n}$,
	\begin{enumerate}
	\item  $n= 3$ or 4, the system \eqref{LaneEmden} has no positive classical solutions.
	
	\item  $n\geq 5$, $pq>1$,  along with
	\begin{align}
		\max\{\alpha,\beta\} > n-3,
		\end{align}
	then system \eqref{LaneEmden} has no positive classical solutions.
	\end{enumerate}
	
	\end{theorem}
For higher order HLS type systems \eqref{pde}, Arthur, Yan and Zhao (2014) \cite{AYZ14} have proved a Liouville theorem with an adapted idea of measure
and feedback argument developed by Souplet (2009) \cite{Souplet09}.

Existence of solutions to subcritical Lane-Emden system with double bounded coefficients can be extended from low dimension to system in higher dimension which fails the integral estimates \eqref{energyEstimates}. This may imply that, integral estimates are essential to problem of nonexistence of solutions. In \cite{CHL14}, Cheng, Huang and Li raised the conjecture below, which is proven to be equivalent to the Lane-Emden conjecture.
\begin{conjecture}
For solution $(u,v)$ to the Lane-Emden system with $p\geq q$, there exist an $s>0$ such that $n-s\beta<1$ and
\begin{align}
	\int_{B_R} v^s \leq C R^{n-s\beta}.  \label{vEstimateBetter}
\end{align}
\end{conjecture}

	\begin{theorem}[Cheng-Huang-Li-2014]\label{CHL14theorem}
		Let $(u,v)$ be a non-negative bounded solution to \eqref{LaneEmden}, with
		$\frac{1}{p+1} +\frac{1}{q+1} > \frac{n-2}{n}$.
		Assume additionally $p\geq q$ and there exists an $s>0$ satisfying $n-s\beta<1$ such that v satisfies \eqref{vEstimateBetter}
%		\begin{align}
%		    \int_{B_R} v^s \leq CR^{n-s\beta}, 
%		\end{align}
		then $u, v \equiv 0$.
	\end{theorem}
	
\begin{remark}
Energy estimate \eqref{vEstimateBetter} is a necessary condition to the Lane-Emden conjecture. One just needs to notice that when $u,v\equiv 0$, \eqref{vEstimateBetter} is obviously satisfied. The key to the proof of Theorem \ref{CHL14theorem} is to show \eqref{vEstimateBetter} is sufficient.
\end{remark}
	
\begin{remark}
The assumption on $v$ is weaker than the corresponding assumption on $u$ due to a comparison principle between $u$ and $v$.

Indeed, we can  replace \eqref{vEstimateBetter} by: for some $r>0$, $n-r\alpha <1$,
\begin{align}
	\int_{B_R} u^r \leq CR^{n-r\alpha}.  \label{uEstimateBetter}
\end{align}
\end{remark}

The proof of Theorem \ref{CHL14theorem} is also based on the feedback argument developed by Souplet, though some basic estimates are adapted as needed. Here a
sketch of the proof of Theorem \ref{CHL14theorem} is presented.
\begin{enumerate}
\item Let $F(R) = \int_{B_R} u^{q+1}.$
By a Rellich-Pohozaev identity of \eqref{LaneEmden} on $B_R$, one can estimate $F(R)$ by quantities on $S^{n-1}$,
\begin{align*}
	F(R) \leq G_1(R) + G_2(R),
\end{align*}
where
\begin{align}
	& G_1(R)=R^n\int_{S^{n-1}}u^{q+1}(R), \label{g1}  \\
	    & G_2(R)=R^n\int_{S^{n-1}}\left( |D_x u(R)|+R^{-1}u(R)\right) \left( |D_x v(R)|+R^{-1}v(R)\right) .  \label{g2}
\end{align}

\item Heuristically, we prove that there exists a sequence $\{R_j\} \rightarrow \infty$ such that the following estimate holds,
\begin{align*}
	F(R_j) \leq C R_j^{-a}F^b(R_j),
\end{align*}
with $a>0$ and $b<1$. Then $F(R)\equiv 0$.
\end{enumerate}

	We start with estimate on $G_1(R)$. By H\"{o}lder's inequality and Sobolev embedding on $S^{n-1}$,
	\begin{align}\label{holderInequalityG1}
	    \|u\|_{q+1} &\leq \|u\|_{\lambda}^{\theta}\|u\|_{\mu}^{1-\theta} \\
	                &\leq (R^2\|D_x^2 u\|_l+\|u\|_1)^{\theta}(R^2\|D_x^2 u\|_k+\|u\|_1)^{1-\theta},
	\end{align}
	So,
	\begin{align*}
		G_1(R) \leq R^n \left( (R^2\|D_x^2 u\|_l+\|u\|_1)^{\theta}(R^2\|D_x^2 u\|_k+\|u\|_1)^{1-\theta}\right) ^{q+1}.
	\end{align*}
	Then by $W^{2,p}$-estimate, energy estimates in Lemma \ref{energyEstimates} together with the assumed integral estimate we have, $\exists \tilde{R}\in(R,2R)$ such that
	\begin{align*}
	    G_1(\tilde{R})\leq R^{-a}F^b(4R).
	\end{align*}
	
	\begin{remark}
	The existence of such $\tilde{R}$  is guaranteed by a fact that for $f\in L_{loc}^p(B_R)$, $\exists \tilde{R}\in(R,2R)$ such that
	\begin{align*}
		\|f\|_{L^{p}(S^{n-1})}(\tilde{R})  \leq M R^{-\frac{n}{p}}\|f\|_{L^{p}(B_{2R})}.
	\end{align*}
	Notice that the same conclusion can be made to finitely many $L^{p_i}_{loc}$ functions.
	\end{remark}

Moreover, we have the following relation between $a$ and $b$,
	\begin{align*}
		  a = (\alpha+\beta+2-n)(1-b).
	\end{align*}
So, to show that $a<0$ and $b<1$ we only need to verify that $b<1$, and this is guaranteed by
	\begin{align*}
%		& b < 1 \\
%		 \Leftrightarrow& (1-\theta)p(q+1) < p+1 \\
%		 \Leftrightarrow& \frac{\frac{1}{l}-\frac{2}{n-1}-\frac{1}{1+q}}{\frac{1}{l}-\frac{1}{k}}(q+1)<k \\
%		 \Leftrightarrow& (\frac{1}{l}-\frac{2}{n-1})(q+1) -1< \frac{k}{l}-1 \\
%		 \Leftrightarrow& \frac{1}{l}(q+1-1-\frac{1}{p}) < \frac{2}{n-1}(q+1) \\
%		 \Leftrightarrow& \frac{pq-1}{s} < \frac{2(q+1)}{n-1} \\
%		 \Leftrightarrow& n-1 < s\beta.
		b<1 \Leftrightarrow n-1 < s\beta.
	\end{align*}
Then we look at estimate on $G_2$.
By H\"{o}lder's inequality, \eqref{g2} becomes,
\begin{align*}
	G_2(R) &\leq R^n\||D_x u| + R^{-1}u\|_z\||D_x v| + R^{-1}v\|_{z'} \\
				&\leq R^n(\|D_x u\|_z+R^{-1}\|u\|_1)(\|D_x v\|_{z'} + R^{-1}\|v\|_1),
\end{align*}
Again applying  H\"{o}lder's inequality and Sobolev embedding to $\|D_x u\|_z$ and $\|D_x v\|_{z'}$, we get
%\begin{align*}
%	\|D_x u\|_z &\leq \|D_x u\|_{\rho_1}^{\tau_1}\|D_x u\|_{\gamma_1}^{1-\tau_1} \\
%						&\leq (R\|D^2_x u\|_l + \|D_x u\|_1)^{\tau_1}(R\|D^2_x u\|_k+ \|D_x u\|_1)^{1+\tau_1}, \\		
%	\|D_x v\|_{z'} &\leq \|D_x v\|_{\rho_2}^{\tau_2}\|D_x v\|_{\gamma_2}^{1-\tau_2} \\
%							&\leq (R\|D^2_x v\|_{1+\epsilon} + \|D_x v\|_1)^{\tau_2}(R\|D^2_x v\|_m+ \|D_x v\|_1)^{1+\tau_2},					
%\end{align*}
%So,
\begin{align*}
    G_2(R)\leq & C R^{n+2}(\|D_x^2 u\|_{l} + R^{-1}\|D_x u\|_1 + R^{-2}\|u\|_1)^{\tau_1} \\
            & \times(\|D_x^2 u\|_k + R^{-1}\|D_x u\|_1 + R^{-2}\|u\|_1)^{1-\tau_1} \\
            & \times(\|D_x^2 v\|_{1+\epsilon} + R^{-1}\|D_x v\|_1 + R^{-2}\|v\|_1)^{\tau_2} \\
            & \times(\|D_x^2 v\|_m + R^{-1}\|D_x v\|_1 + R^{-2}\|v\|_1)^{1-\tau_2},
\end{align*}
with
\begin{align}\label{rangeZ}
	\max\left\lbrace \frac{1}{k} - \frac{1}{n-1}, \frac{1}{n-1}\right\rbrace \leq \frac{1}{z} \leq \min \left\lbrace \frac{1}{l}-\frac{1}{n-1}, \frac{1}{q+1}+\frac{1}{n-1}\right\rbrace.
\end{align}
Similar to $G_1(R)$,
by $W^{2,p}$-estimate, basic integral estimates \eqref{uvEstimate1} and \eqref{uvEsitmatePQ} together with the assumed integral estimate we have, $\exists \tilde{R}\in(R,2R)$ such that
	\begin{align*}
	    G_2(\tilde{R})\leq R^{-\tilde{a}}F^{\tilde{b}}(4R).
	\end{align*}
(Surprisingly) similar to $a,b$ for $G_1(R)$, $\tilde{a}$ and $\tilde{b}$ also have the relation
	\begin{align*}
		\tilde{a} = (\alpha+\beta+2-n)(1-\tilde{b}),
	\end{align*}
	so $\tilde{a}>0$ if and only if $\tilde{b}<1$,
	and
	\begin{align}
	\tilde{b} &< 1 \label{bLessThan1} \\
	\Leftrightarrow & (m-\frac{k}{l})\frac{1}{z}+(\frac{k}{n-1}+(m-1)(k-1))\frac{1}{l}+\frac{m-2}{n-1}-(m-1)>0, \label{zAndl}
	\end{align}
	Now we only need to verify this inequality holds when $\frac{1}{z}$ takes value as its upper bounds in \eqref{rangeZ}:
	\begin{align*}
		\frac{1}{z} \leq \min \left\lbrace \frac{1}{l}-\frac{1}{n-1}, \frac{1}{q+1}+\frac{1}{n-1}\right\rbrace.
	\end{align*}

\textbf{Case 1.} If $\frac{1}{l}-\frac{1}{n-1} \geq \frac{1}{q+1}+\frac{1}{n-1}$,
\begin{align*}
 \tilde{b} < 1
 	\Leftrightarrow  n- s\beta  < 1,
\end{align*}
which is true under our assumption.

\textbf{Case 2.} If $\frac{1}{l}-\frac{1}{n-1} < \frac{1}{q+1}+\frac{1}{n-1}$,
\begin{align*}
 \tilde{b} < 1
 	\Leftrightarrow  \frac{1}{k} < \frac{1}{l} < \frac{2}{n-1}+\frac{1}{q},
\end{align*}
where the inequalities on the right are true since $\frac{1}{l} < \frac{1}{q+1}+\frac{2}{n-1}< \frac{2}{n-1}+\frac{1}{q}$.

 In all, \eqref{bLessThan1} always holds under our assumption $n-s\beta<1$. This finishes the proof of Theorem \ref{CHL14theorem}.

\newpage

\providecommand{\bysame}{\leavevmode\hbox to3em{\hrulefill}\thinspace}
\providecommand{\MR}{\relax\ifhmode\unskip\space\fi MR }
\providecommand{\MRhref}[2]{%
  \href{http://www.ams.org/mathscinet-getitem?mr=#1}{#2}
}
\providecommand{\href}[2]{#2}

\noindent\href{mailto:ze.cheng@colorado.edu}{ze.cheng@colorado.edu}\\
\href{mailto:genggenghuang1986@gmail.com}{genggenghuang1986@gmail.com}\\
\href{mailto:congmingli@gmail.com}{congmingli@gmail.com}
\end{document}